\documentclass{article}
\usepackage{amssymb}
\usepackage{amsmath}
\usepackage{amsfonts}
\usepackage{graphicx}%
\providecommand{\U}[1]{\protect \rule{.1in}{.1in}}

\newtheorem{theorem}{Theorem}

\newtheorem{corollary}[theorem]{Corollary}

\newtheorem{proposition}[theorem]{Proposition}
\newtheorem{remark}[theorem]{Remark}

\begin{document}

\title{\textbf{On the Tchebychev Vector Field}\\ \textbf{in the Relative Differential Geometry }}
\author{\textbf{Stylianos Stamatakis\strut \medskip}\\ \emph{Aristotle University of Thessaloniki}\\ \emph{Department of Mathematics}\\ \emph{GR-54124 Thessaloniki, Greece}\\ \emph{e-mail: stamata@math.auth.gr}}
\date{}
\maketitle

\begin{abstract}
\noindent In this paper we deal with relative normalizations of
hypersurfaces\ in the Euclidean space $%
\mathbb{R}
^{n+1}.$ Considering a relative normalization $\bar{y}$ of an hypersurface
$\Phi$ we decompose the corresponding Tchebychev vector $\bar{T}$ in two
components, one parallel to the Tchebychev vector $\bar{T}_{EUK}$ of the
Euclidean normalization $\bar{\xi}$ and one parallel to the orthogonal
projection $\bar{y}_{T}$ of $\bar{y}$ in the tangent hyperplane of $\Phi$. We
use this decomposition to investigate some properties of $\Phi$, which concern
its Gaussian curvature, the support function, the Tchebychev vector field
etc.\smallskip \newline MSC 2010: 53A05, 53A07, 53A15, 53A40 \smallskip \newline Keywords:
Euclidean hypersurfaces, Relative normalizations, equiaffine normalization, Tchebychev vector

\end{abstract}

\section{Introduction}

\noindent To set the stage for this work the classical notation of relative
differential geometry is briefly presented; for this purpose the paper
\cite{Schneider2} is used as general reference.\smallskip

\noindent In the Euclidean space $%
\mathbb{R}
^{n+1}$ let $\Phi:\bar{x}=\bar{x}(u^{i}),(u^{i}):=(u^{1},\cdots,u^{n})\in
U\subset%
\mathbb{R}
^{n}$ be an injective $C^{r}$-immersion with Gaussian curvature $K\neq0$
$\forall$ $(u^{i})\in U$. A $C^{s}$-mapping $\bar{y}:U\longrightarrow%
\mathbb{R}
^{n+1}$ ($r>s\geq1$) is called a \textit{relative }$C^{s}$%
-\textit{normalization} if%
\begin{equation}
\bar{y}(P)\notin T_{P}\Phi,\quad \bar{y}_{/i}(P)\in T_{P}\Phi \quad
(i=1,\ldots,n)~~\footnote{Partial derivatives of a function $f$ are denoted by
$f_{/i}:=\frac{\partial f}{\partial u^{i}},~f_{/ij}:=\frac{\partial^{2}%
f}{\partial u^{i}\partial u^{j}}\medskip$} \label{1}%
\end{equation}
at every point $P\in \Phi$, where $T_{P}\Phi$ is the tangent vector space of
$\Phi$ at $P$.\smallskip

\noindent The \textit{covector }$\bar{X}$ of the tangent hyperplane is defined
by%
\begin{equation}
\langle \bar{X},\bar{x}_{/i}\rangle=0\quad(i=1,\ldots,n)\quad \text{and}%
\quad \langle \bar{X},\bar{y}\rangle=1, \label{2}%
\end{equation}
where $\langle~,~\rangle$\ denotes the standard scalar product in $%
\mathbb{R}
^{n+1}$. The \textit{relative metric} $G$ on $U$ is introduced by
\begin{equation}
G_{ij}=\langle \bar{X},\bar{x}_{/ij}\rangle. \label{4}%
\end{equation}
The \textit{support function }of the relative normalization $\bar{y}$ is
defined by%
\[
q:=\langle \bar{\xi},\bar{y}\rangle:U\longrightarrow%
\mathbb{R}
,\quad q\in C^{s}\left(  U\right)  ,
\]
where $\bar{\xi}:U\longrightarrow%
\mathbb{R}
^{n+1}$ is the \textit{Euclidean normalization} of $\Phi$. Because of
(\ref{1}) the support function $q$ never vanishes on $U$. Furthermore, from
(\ref{2}) it follows the relation
\begin{equation}
\bar{X}=q^{-1}\bar{\xi}. \label{6}%
\end{equation}
On account of (\ref{4}) and (\ref{6}), we obtain%
\begin{equation}
G_{ij}=q^{-1}h_{ij}, \label{7}%
\end{equation}
where $h_{ij}$ are the components of the second fundamental form $II$ of
$\Phi$. We mention that given a support function $q,$ the relative
normalization $\bar{y}$ is uniquely determined and possesses the following
parametrization (see \cite[p. 197]{Manhart})
\begin{equation}
\bar{y}=-h^{ij}\,q_{/i}\, \bar{x}_{/j}+q\, \bar{\xi}, \label{8}%
\end{equation}
where $h^{ij}$ are the components of the inverse tensor of $h_{ij}$
\footnote{From now on we follow similar notation for the inverse of a given
tensor}.\smallskip

\noindent Let $Q$ be a definite quadratic form. For a $C^{r}(U)$-function $f$
we denote by $\nabla^{Q}f$ the first Beltrami differential operator, by
$\triangle \!^{Q}f$ the second Beltrami differential operator and by
$\nabla_{i}^{Q}f$ the covariant derivative of $f$ with respect to
$Q$.\smallskip

\noindent We consider the components
\[
A_{ijk}:=\langle \bar{X},\nabla_{k}^{G}\nabla_{j}^{G}\bar{x}_{/i}\rangle
\]
of the symmetric \textit{Darboux tensor}. Then the \textit{Tchebychev vector
field }$\bar{T}$, which corresponds to the\thinspace relative normalization
$\bar{y},$ is defined by%
\[
\bar{T}=\frac{1}{n}A_{ijk}\,G^{jk}\,G^{im}\, \bar{x}_{/m}.
\]
We mention that the relation
\begin{equation}
\frac{\triangle^{G}\bar{x}}{n}=\bar{T}+\bar{y}. \label{11}%
\end{equation}
holds \cite{Heil}. The \textit{relative shape operator} $B$ has the components
$B_{i}^{j}:U\longrightarrow%
\mathbb{R}
,$ defined by%
\[
\bar{y}_{/i}=-B_{i}^{j}\, \bar{x}_{/j}.
\]
For the \textit{relative mean curvature} $H:=\frac{1}{n}B_{i}^{i}$ we have
according to \cite{Manhart}%
\begin{equation}
H=q\,H_{I}+\frac{1}{n}\left[  \triangle \!^{II}\left(  \ln q\right)
+\frac{2-n}{4}\triangledown^{II}\left(  \ln q\right)  \right]  , \label{41}%
\end{equation}
where $H_{I}$ is the Euclidean mean curvature of $\Phi$.

\section{A decomposition of the Tchebychev vector field}

\noindent Let $I=g_{ij}du^{i}du^{j}$ be the first fundamental form and
$III=e_{ij}du^{i}du^{j}$ the third fundamental form of $\Phi$. Taking into
account the Weingarten equations%
\[
\bar{\xi}_{/i}=-h_{ij}\,g^{jk}\, \bar{x}_{/k},
\]
and the relations (\ref{7}), we obtain%
\begin{equation}
\nabla^{II}(f,\bar{\xi})=-\nabla^{I}(f,\bar{x})=q^{-1}\, \nabla^{G}(f,\bar
{\xi}), \label{16b}%
\end{equation}%
\begin{equation}
\nabla^{II}(f,\bar{x})=-\nabla^{III}(f,\bar{\xi})=q^{-1}\, \nabla^{G}%
(f,\bar{x}). \label{16a}%
\end{equation}
\noindent We firstly compute the vectors $\triangle \!^{II}\bar{x}$ and
$\triangle \!^{II}\bar{\xi}$. To this end we consider the components
$T_{jk}^{i}={}^{I}\Gamma_{jk}^{i}-{}^{II}\Gamma_{jk}^{i}$ of the
\textit{difference tensor }of the Levi-Civita connections with respect to $I$
and $II$. It is known that the relations%
\[
T_{ij}^{k}=-\frac{1}{2}h^{km}\, \nabla_{m}^{I}h_{ij},
\]%
\[
\nabla_{j}^{II}\bar{x}_{/i}=T_{ij}^{k}\, \bar{x}_{/k}+h_{ij}\, \bar{\xi
}\mathbf{,}%
\]%
\[
T_{im}^{i}=\frac{-K_{/m}}{2K}%
\]
hold on $U$ (see \cite[p. 22]{Huck} and \cite[p.197]{Manhart}). Using them and
the Mainardi-Codazzi equations $\nabla_{m}^{I}h_{ij}-\nabla_{j}^{I}h_{im}=0$
we find%
\[
-2h^{ij}\,T_{ij}^{k}=h^{ij}\,h^{km}\, \nabla_{m}^{I}h_{ij}=h^{km}\,h^{ij}\,
\nabla_{j}^{I}h_{im}=-2h^{km}\,T_{im}^{i}=\frac{h^{km}}{K}\,K_{/m},
\]
so that, by a direct computation, we arrive at
\begin{equation}
\triangle \!^{II}\bar{x}=\frac{-1}{2K}\nabla^{II}(K,\bar{x})+n\, \bar{\xi}.
\label{13}%
\end{equation}

\noindent Following similar computation and taking account of the relations
(for $n=2$ see \cite{Huck}),%
\[
\nabla_{j}^{II}\bar{\xi}_{/i}=-T_{ij}^{k}\, \bar{\xi}_{/k}-e_{ij}\, \bar{\xi
},
\]%
\[
H_{I}=\frac{e_{ij}\,h^{ij}}{n},
\]
we find%
\begin{equation}
\triangle \!^{II}\bar{\xi}=\frac{1}{2K}\nabla^{II}(K,\bar{\xi})-n\,H_{I}\,
\bar{\xi}. \label{14}%
\end{equation}

\begin{remark}
Relations (\ref{13}) and (\ref{14}) for the case $n=2$ were proved in
\cite{Stamatakis} with sign convention such that $\triangle=-\frac
{\partial^{2}}{\partial x^{2}}-\frac{\partial^{2}}{\partial y^{2}}$ for the
metric $ds^{2}=dx^{2}+dy^{2}$.
\end{remark}

\noindent Continuing our considerations, we may compute the second Beltrami
differential operator $\triangle^{G}f$ for a $C^{r}(U)$-function $f$. It is
known (\cite[p. 196]{Manhart}) that between the Levi-Civita connections with
respect to $G$ and $II$ the relation%
\[
^{G}\! \Gamma_{ij}^{k}=\!^{II}\! \Gamma_{ij}^{k}-\frac{\delta_{i}^{k}%
\,q_{/j}+\delta_{j}^{k}\,q_{/i}-h_{ij}\,h^{km}\,q_{/m}}{2q}%
\]
holds. By using them and (\ref{7}), we obtain the relation%
\begin{equation}
\triangle \!^{G}f=q\, \triangle^{II}f-\frac{n-2}{2}\nabla^{II}(q,f). \label{15}%
\end{equation}
We apply now (\ref{15}) to $\bar{x}$ and make use of (\ref{8}), (\ref{16a})
and (\ref{13}). Thus we find%
\begin{equation}
\triangle \!^{G}\bar{x}=\frac{q}{2K}\nabla^{III}(K,\bar{\xi})-\frac{n+2}%
{2}\nabla^{III}(q,\bar{\xi})+n\, \bar{y}. \label{18}%
\end{equation}
Similarly, by applying (\ref{15}) to $\bar{\xi}$ and by using (\ref{16b}) and
(\ref{14}) we get%
\begin{equation}
\triangle \!^{G}\bar{\xi}=\frac{-q}{2K}\nabla^{I}(K,\bar{x})+\frac{n-2}%
{2}\nabla^{I}(q,\bar{x})-n\,q\,H_{I}\, \bar{\xi}. \label{19}%
\end{equation}
Taking into account (\ref{11}) and (\ref{18}) we obtain for the Tchebychev
vector field $\bar{T}$ of $\Phi$, which corresponds to the support function
$q,$
\begin{equation}
\bar{T}=\frac{q}{2nK}\nabla^{III}(K,\bar{\xi})-\frac{n+2}{2n}\nabla
^{III}(q,\bar{\xi}). \label{19b}%
\end{equation}
In the case of the Euclidean normalization ($q=1$) it is $\bar{y}=\bar{\xi}$,
whereupon we find for the corresponding Tchebychev vector field
\begin{equation}
\bar{T}_{{\small _{EUK}}}=\frac{1}{2nK}\nabla^{III}(K,\bar{\xi}). \label{20}%
\end{equation}
Introducing the tangent vector field%
\begin{equation}
\bar{Q}:=\frac{1}{2nq}\nabla^{III}(q,\bar{\xi}) \label{21}%
\end{equation}
of $\Phi$ and inserting this, as well as (\ref{20}), in (\ref{19b}), we get%
\begin{equation}
\bar{T}=q\, \bar{T}_{{\small _{EUK}}}-q\,(n+2)\, \bar{Q}. \label{22}%
\end{equation}
Similarly, we obtain from (\ref{18})%
\begin{equation}
\triangle \!^{G}\bar{x}=n\,q\,[\bar{T}_{{\small _{EUK}}}+(n-2)\, \bar{Q}%
+\bar{\xi}] \label{24}%
\end{equation}
and from (\ref{8}), (\ref{16a})%
\begin{equation}
\bar{y}=q\,(2n\, \bar{Q}+\bar{\xi}), \label{23}%
\end{equation}
i.e. the vector $\bar{Q}$ is parallel to the orthogonal projection $\bar
{y}_{{_{T}}}$ of the relative normalization $\bar{y}$ in the the tangent
vector space $T_{P}\Phi$ of $\Phi$ at $P$.\smallskip

\noindent From (\ref{22}) and (\ref{23}) we see that \textit{the vector field
}$\frac{1}{q}\left[  2n\, \bar{T}+\left(  n+2\right)  \, \bar{y}\right]  $
\textit{is independent of the \thinspace relative normalization and equals}
$2n$\thinspace$\bar{T}_{{\small _{EUK}}}+\left(  n+2\right)  $\thinspace
$\bar{\xi}$.$\smallskip$

\noindent Finally, taking into account (\ref{16a}) and (\ref{19b}), we can
write the Tchebychev vector field $\bar{T}$ as gradient (see \cite[p.
243]{Heil})%
\begin{equation}
\bar{T}=\nabla \!^{G}(\ln \varphi(u^{i}),\bar{x}), \label{25}%
\end{equation}
where
\begin{equation}
\varphi(u^{i})=|K|^{\frac{-1}{2n}}\cdot|q|^{\frac{n+2}{2n}}. \label{26}%
\end{equation}
Consequently, $\bar{T}$ \textit{is irrotational with respect to the relative
metric }$G$.

\section{Relatively normalized surfaces by $^{\left(  \alpha \right)  }\bar{y}%
$}

\noindent We consider now the relative normalizations $^{\left(
\alpha \right)  }\bar{y}:U\longrightarrow%
\mathbb{R}
^{n}$, which are introduced by F. Manhart \cite{Manhart}, and, on account of
(\ref{8}), are defined by the support functions%
\begin{equation}
^{_{(\alpha)}}q:=|K|^{\alpha},\quad \alpha \in%
\mathbb{R}
. \label{27}%
\end{equation}
Denoting by $^{_{(\alpha)}}\bar{Q}$ the corresponding vector field and taking
into account (\ref{20}) and (\ref{21}), we find%
\begin{equation}
^{_{(\alpha)}}\bar{Q}=\alpha \, \bar{T}_{_{EUK}}. \label{42}%
\end{equation}
Conversely, if the vector fields $\bar{Q}$ and $\bar{T}_{{\small _{EUK}}}$ are
such that $\bar{Q}=\alpha$\thinspace$\bar{T}_{{\small _{EUK}}}$ for $\alpha \in%
\mathbb{R}
,$ it turns out that
\begin{equation}
q=\lambda \,|K|^{\alpha}, \label{33}%
\end{equation}
where $\lambda \in%
\mathbb{R}
^{\ast}:=%
\mathbb{R}
-\left \{  0\right \}  $ is an arbitrary constant. Consequently, we have the following

\begin{proposition}
The vector fields $\bar{Q}$ and $\bar{T}_{{\small _{EUK}}}$ satisfy the
relation $\bar{Q}=\alpha$\thinspace$\bar{T}_{{\small _{EUK}}}$ for a constant
$\alpha \in%
\mathbb{R}
,$ if and only if the support function \textit{has the form }$q=\lambda
$\thinspace$|K|^{\alpha},$ where $\lambda \in%
\mathbb{R}
^{\ast}$ is an arbitrary constant.
\end{proposition}

\noindent We denote by $^{_{(\alpha)}}G\ $the relative metric and by
$^{_{(\alpha)}}\bar{T}$ the Tchebychev vector with respect to the relative
normalization $^{_{(\alpha)}}\bar{y}.$ Then from (\ref{22}), (\ref{24}),
(\ref{23}) and (\ref{42}) we have%
\[
^{_{(\alpha)}}\bar{y}=\ ^{_{(\alpha)}}\!q\,[2\alpha n\, \bar{T}_{_{EUK}}%
+\bar{\xi}],
\]%
\begin{equation}
^{_{(\alpha)}}\bar{T}=\ ^{_{(\alpha)}}\!q\,[1-\alpha(n+2)]\, \bar{T}_{_{EUK}},
\label{30}%
\end{equation}%
\begin{equation}
\triangle^{^{_{(\alpha)}}G}\bar{x}=n\ ^{_{(\alpha)}}\!q\, \left \{
[1+\alpha \,(n-2)]\, \bar{T}_{_{EUK}}+\bar{\xi}\right \}  , \label{31}%
\end{equation}
while for the function $\varphi(u^{i})$ in (\ref{26}) we find%
\[
^{_{(\alpha)}}\varphi(u^{i})=|K|^{\frac{\alpha(n+2)-1}{2n}}.
\]
We note that the formulae (\ref{30})-(\ref{31}) remain invariant if the
support function $q$ has the form (\ref{33}).\smallskip

\noindent In the one-parameter family of relative normalizations $^{\left(
\alpha \right)  }\bar{y},$ which are determined by the support functions
$^{\left(  \alpha \right)  }q$, among other relative normalizations

\begin{itemize}
\item the \textit{Euclidean normalization }(when $\alpha=0$) and \smallskip

\item the \textit{equiaffine normalization} (when $\alpha=1/\left(
n+2\right)  $)
\end{itemize}

\noindent are contained. Furthermore we find%
\[
^{_{(0)}}\bar{Q}=\bar{0},\quad^{_{(0)}}\bar{y}=\bar{\xi},\quad^{_{(0)}}\bar
{T}=\bar{T}_{_{EUK}},
\]%
\[
^{(\frac{1}{n+2})}\bar{Q}=\frac{1}{n+2}\bar{T}_{_{EUK}},\quad^{(\frac{1}%
{n+2})}\bar{y}=\sqrt[n+2]{|K|}\ [\frac{2n}{n+2}\bar{T}_{_{EUK}}+\bar{\xi
}],\quad^{(\frac{1}{n+2})}\bar{T}=\bar{0}.
\]

\section{Applications}

\textbf{4.1. }In this paragraph, using the vector fields $\bar{T}%
,\triangle^{G}\bar{x}$ and $\triangle^{G}\bar{\xi},$ we find necessary and
sufficient conditions for the Gaussian curvature $K$ of $\Phi$ to be constant
or for the support function $q$ to be of the form (\ref{33}).\bigskip

\noindent \textbf{A. }From (\ref{20}) and (\ref{21}) we have
\[
\langle \bar{T}_{_{EUK}},d\bar{\xi}\rangle=\frac{1}{2nK}dK,
\]%
\[
\langle \bar{Q},d\bar{\xi}\rangle=\frac{1}{2nq}dq,
\]
so that, from (\ref{24}) we obtain%
\begin{equation}
2\langle \triangle \!^{G}\bar{x},d\bar{\xi}\rangle=q\,d\left(  \ln \left(
|K|\cdot|q|^{n-2}\right)  \right)  . \label{35}%
\end{equation}
Hence, we have

\begin{proposition}
(a) \textit{When }$n=2,$ it holds $\langle \triangle \!^{G}\bar{x},d\bar{\xi}\rangle
=0$\textit{ }at every point $P\in \Phi$\textit{ if and only if }$K=const$%
\textit{.}\smallskip \newline(b) \textit{When }$n\geq3,$ the relation\textit{
}$\langle \triangle \!^{G}\bar{x},d\bar{\xi}\rangle=0$\textit{ holds }at every
point $P\in \Phi$\textit{ if and only if the support function has the form
}$q=\lambda$\thinspace$|K|^{\frac{1}{2-n}},$ where $\lambda$ is an arbitrary
not vanishing constant\textit{.\smallskip}
\end{proposition}

\noindent \textbf{B. }From (\ref{16a}) and (\ref{25}) it follows%
\[
\langle \bar{T},d\bar{\xi}\rangle=-q\,d(\ln \varphi).
\]
Obviously, we have $\langle \bar{T},d\bar{\xi}\rangle=0$ if and only if the
function $\varphi(u^{i}),$ which is defined in (\ref{26}), is constant, or if
and only if
$\vert$%
$K|\cdot|q|^{-(n+2)}=const$., i.e. the support function has the form%
\[
q=\lambda \,|K|^{\frac{1}{n+2}}%
\]
where $\lambda \in%
\mathbb{R}
^{\ast}$ is an arbitrary constant. So we have

\begin{proposition}
The following properties are equivalent:\smallskip \newline(a) The function
$\varphi(u^{i})=|K|^{\frac{-1}{2n}}\cdot|q|^{\frac{n+2}{2n}}$ is
constant.\smallskip \newline(b) $\langle \bar{T},d\bar{\xi}\rangle=0\ $at every
point $P\in \Phi.$\smallskip \newline(c) The support function has the form
$q=\lambda$\thinspace$|K|^{\frac{1}{n+2}},$ where $\lambda \in%
\mathbb{R}
^{\ast}$ is an arbitrary constant.\smallskip
\end{proposition}

\noindent \textbf{C. }From (\ref{19}) we have%
\begin{equation}
2\langle \triangle \!^{G}\bar{\xi},d\bar{x}\rangle=q\,d(\ln(|K|^{-1}%
\cdot|q|^{n-2})). \label{36}%
\end{equation}
Consequently, we obtain

\begin{proposition}
(a) When $n=2,$ it holds $\langle \triangle^{G}\bar{\xi},d\bar{x}\rangle=0$ at
every point $P\in \Phi$ if and only if $K=const$.\smallskip \newline(b) When
$n\geq3,$ the relation $\langle \triangle \!^{G}\bar{\xi},d\bar{x}\rangle=0$
holds at every point $P\in \Phi$ if and only if the support function has the
form $q=\lambda$\thinspace$|K|^{\frac{1}{n-2}},$ where $\lambda \in%
\mathbb{R}
^{\ast}$ is an arbitrary constant.\smallskip
\end{proposition}

\noindent \textbf{D. }From (\ref{35}) and (\ref{36}) it follows%
\[
\langle \triangle \!^{G}\bar{x},d\bar{\xi}\rangle+\langle \triangle \!^{G}\bar
{\xi},d\bar{x}\rangle=(n-2)\,dq,
\]%
\[
\langle \triangle \!^{G}\bar{x},d\bar{\xi}\rangle-\langle \triangle \!^{G}\bar
{\xi},d\bar{x}\rangle=\frac{q}{K}dK,
\]
which lead to

\begin{proposition}
(a) Let $n=2.$ For every relative normalization the relation $\langle
\triangle \!^{G}\bar{x},d\bar{\xi}\rangle+\langle \triangle \!^{G}\bar{\xi}%
,d\bar{x}\rangle=$ $0$ holds at every point $P\in \Phi$.\smallskip \newline(b)
For $n\neq2$ the relation $\langle \triangle \!^{G}\bar{x},d\bar{\xi}%
\rangle+\langle \triangle \!^{G}\bar{\xi},d\bar{x}\rangle=$ $0$ holds at every
point $P\in \Phi$ if and only if the support function $q$ is
constant.\smallskip \newline(c) The relation $\langle \triangle \!^{G}\bar
{x},d\bar{\xi}\rangle-\langle \triangle \!^{G}\bar{\xi},d\bar{x}\rangle=$ $0$
holds at every point $P\in \Phi$ if and only if the Gaussian curvature $K$ is
constant.\smallskip
\end{proposition}

\noindent \textbf{4.2. }Let $\bar{y}_{i},i=1,2,$ be two relative normalizations
of $\Phi$. We denote by $q_{i},G_{i}$ and $\bar{T}_{i}$ the corresponding
support functions, relative metrics and Tchebychev vector fields,
respectively.\bigskip

\noindent \textbf{A. }From (\ref{23}) it turns out%
\[
\langle \bar{y}_{1},\bar{y}_{2}\rangle=\nabla^{III}(q_{1},q_{2})+q_{1}\,q_{2},
\]
so that: \textit{The \thinspace relative normalizations }$\bar{y}_{1}$
\textit{and} $\bar{y}_{2}$ \textit{are orthogonal if and only if the
corresponding support functions satisfy the relation}%
\[
\nabla^{III}(\ln|q_{1}|,\ln|q_{2}|)=-1.\medskip
\]
\textbf{B. }On account of (\ref{21}), we have:

\begin{proposition}
The \textit{vector fields }$\bar{Q}_{1}$\textit{ and }$\bar{Q}_{2}$\textit{
satisfy the relation }$\bar{Q}_{2}=\alpha$\thinspace$\bar{Q}_{1}$
\textit{for\ a constant }$\alpha \in%
\mathbb{R}
$, if and only if \textit{the corresponding support functions }$q_{1}$\textit{
and }$q_{2}$\textit{ satisfy the relation }$q_{2}=\lambda$\thinspace
$q_{1}^{\alpha},$ where $\lambda \in%
\mathbb{R}
^{\ast}$ is an arbitrary constant\textit{.}
\end{proposition}

\noindent For the corresponding relative normalizations and Tchebychev vector
fields we find%
\[
\bar{y}_{2}=\lambda \,q_{1}^{\alpha-1}\, \left[  \bar{y}_{1}+\left(
\alpha-1\right)  \, \nabla^{III}(q_{1},\bar{\xi})\right]  ,
\]%
\[
\bar{T}_{2}=\lambda \,q_{1}^{\alpha-1}\left[  \bar{T}_{1}-\frac{\left(
\alpha-1\right)  \left(  n+2\right)  }{2n}\nabla^{III}(q_{1},\bar{\xi
})\right]  .\medskip
\]
\noindent \textbf{C. }We study now the case of two relative normalizations
$\bar{y}_{i},i=1,2,$ for which there is a constant $\alpha \in%
\mathbb{R}
$, such that the corresponding Tchebychev vector fields satisfy the relation
$\bar{T}_{2}=\alpha$\thinspace$\bar{T}_{1}$. Taking into account (\ref{22}),
we see that the last relation is equivalent to%
\begin{equation}
(q_{2}-\alpha \,q_{1})\, \nabla^{III}(\ln|K|,\bar{\xi})-\left(  n+2\right)  \,
\nabla^{III}(q_{2}-\alpha \,q_{1},\bar{\xi})=\bar{0}. \label{38}%
\end{equation}
For $q_{2}\neq \alpha$\thinspace$q_{1}$ it follows that%
\begin{equation}
|q_{2}-\alpha \,q_{1}|=\lambda \,q_{_{AFF}}, \label{39}%
\end{equation}
where $q_{_{AFF}}:=|K|^{\frac{1}{n+2}}$ is the support function of the
equiaffine normalization and $\lambda$ is an arbitrary
positive\ constant\textit{.} For $q_{2}=\alpha$\thinspace$q_{1}$ the relation
(\ref{39})\ is still valid (for $\lambda=0)$. \smallskip

\noindent We denote by $\bar{y}_{_{AFF}}$ the equiaffine relative
normalization. Then, by using (\ref{21}) and (\ref{23}), it turns out that
(\ref{39}) holds, if and only if the relative normalizations $\bar{y}_{i}$,
$i=1,2,$ satisfy the relation%
\[
\bar{y}_{2}=\alpha \, \bar{y}_{1}+\mu \, \bar{y}_{_{AFF}},
\]
for $\mu=\varepsilon$\thinspace$\lambda,$ where $\varepsilon
=\operatorname*{sign}\left(  q_{2}-\alpha \,q_{1}\right)  .$ So we have the result:

\begin{proposition}
\label{prop1}The following properties are equivalent:\smallskip \newline(a) The
Tchebychev vector fields $\bar{T}_{1}$ and $\bar{T}_{2}$ satisfy the relation
$\bar{T}_{2}=\alpha$\thinspace$\bar{T}_{1},$ where $\alpha \in%
\mathbb{R}
.$\smallskip \newline(b) The support functions $q_{1}$\textit{ and }$q_{2}%
$\textit{ satisfy the relation }$|q_{2}-\alpha$\thinspace$q_{1}|=\lambda
$\thinspace$q_{_{AFF}},$ where $\lambda$ is an arbitrary
non-negative\ constant.\smallskip \newline(c) \textit{The }relative
normalizations $\bar{y}_{1}$ and $\bar{y}_{2}$ \textit{satisfy the relation
}$\bar{y}_{2}=\alpha$\thinspace$\bar{y}_{1}+\mu$\thinspace$\bar{y}_{_{AFF}},$
where $\mu$ is an arbitrary constant.
\end{proposition}

\begin{remark}
Given a relative normalization $\bar{y}_{1}$, a $C^{s}$-mapping $\bar{y}%
_{2}:U\longrightarrow%
\mathbb{R}
^{n+1}$ satisfying the relation $\bar{y}_{2}=\alpha$\thinspace$\bar{y}_{1}%
+\mu$\thinspace$\bar{y}_{_{AFF}},$ where $\alpha,\mu \in%
\mathbb{R}
,$ is a relative normalization if and only if $\alpha$\thinspace$q_{1}+\mu
$\thinspace${q}_{_{AFF}}\neq0.$
\end{remark}

\noindent From Proposition \ref{prop1} we obtain the

\begin{corollary}
(a) Let $\bar{y}_{1}$ be a relative normalization of $\Phi$ and $\bar{T}_{1}$
be the corresponding Tchebychev vector field. Then all relative normalizations
of the one-parameter family%
\[
\left \{  \bar{y}~/~\bar{y}=\bar{y}_{1}+\mu \, \bar{y}_{_{AFF}},\, \mu \in%
\mathbb{R}
,~\mu \neq-q_{1}\cdot q_{_{AFF}}^{-1}\right \}  ,
\]
have $\bar{T}_{1}$ as common corresponding Tchebychev vector field.\smallskip
\newline(b) All relative normalizations of the one-parameter family%
\[
\left \{  \bar{y}~/~\bar{y}=\bar{\xi}+\mu \, \bar{y}_{_{AFF}},~\mu \in%
\mathbb{R}
,~\mu \neq- q_{_{AFF}}^{-1}\right \}  ,
\]
have $\bar{T}_{_{EUK}}$ as common corresponding Tchebychev vector field.
\end{corollary}

\noindent As immediate consequences of (\ref{11}), (\ref{41}) it follows:

\begin{proposition}
If the Tchebychev vector fields $\bar{T}_{1}$ and $\bar{T}_{2}$ satisfy the
relation $\bar{T}_{2}=\alpha$\thinspace$\bar{T}_{1}$ for $\alpha \in%
\mathbb{R}
,$ then we have\smallskip \newline(a) The vector fields $\triangle \!^{G_{i}%
}\bar{x},$ $i=1,2,$ satisfy the relation
\[
\triangle \!^{G_{2}}\bar{x}=\alpha \, \triangle \!^{G_{1}}\bar{x}+\frac{\mu}%
{n}\bar{y}_{_{AFF}},
\]
where $\mu$ is an arbitrary constant.\smallskip \newline(b) T\textit{he
corresponding }relative mean curvatures $H_{i}$ of the normalizations $\bar
{y}_{i},$ $i=1,2,$ \textit{satisfy the relation }%
\[
H_{2}=\alpha \,H_{1}+\mu \,H_{_{AFF}},
\]
where $H_{_{AFF}}$ is the equiaffine mean curvature and $\mu$ is an arbitrary
constant.\medskip
\end{proposition}

\noindent \textbf{D. }Finally, from (\ref{19}) and (\ref{24}) we have
\[
\frac{\triangle \!^{G_{1}}\bar{x}}{q_{1}}-\frac{\triangle \!^{G_{2}}\bar{x}%
}{q_{2}}=\frac{n-2}{2}\nabla^{III}\left(  \ln \left \vert \frac{q_{1}}{q_{2}%
}\right \vert ,\bar{\xi}\right)  ,
\]%
\[
\frac{\bar{T}_{1}}{q_{1}}-\frac{\bar{T}_{2}}{q_{2}}=-\left(  n+2\right)  \,
\left(  \bar{Q}_{1}-\bar{Q}_{2}\right)  =-\frac{n+2}{2}\nabla^{III}\left(
\ln \left \vert \frac{q_{1}}{q_{2}}\right \vert ,\bar{\xi}\right)  ,
\]%
\[
\frac{\triangle \!^{G_{1}}\bar{\xi}}{q_{1}}-\frac{\triangle \!^{G_{2}}\bar{\xi}%
}{q_{2}}=\frac{n-2}{2}\nabla^{I}\left(  \ln \left \vert \frac{q_{1}}{q_{2}%
}\right \vert ,\bar{x}\right)  ,
\]
so that we conclude

\begin{proposition}
(a) When $n=2,$ the vector fields $q^{-1}\triangle^{G}\bar{x}$ and
$q^{-1}\triangle^{G}\bar{\xi}$ are independent of the normalization.\smallskip
\newline(b) When $n\geq3$ the following properties are equivalent:\smallskip
\newline(i) $\  \qquad \langle \frac{\triangle^{\!G_{1}}\bar{x}}{q_{1}}%
-\frac{\triangle \!^{G_{2}}\bar{x}}{q_{2}},d\bar{\xi}\rangle=0\ $at every point
$P\in \Phi.$\medskip \newline(ii) \qquad$\langle \frac{\bar{T}_{1}}{q_{1}}%
-\frac{\bar{T}_{2}}{q_{2}},d\bar{\xi}\rangle=0$ at every point $P\in \Phi
.$\medskip \newline(iii)\qquad$\langle \frac{\triangle \!^{G_{1}}\bar{\xi}}%
{q_{1}}-\frac{\triangle \!^{G_{2}}\bar{\xi}}{q_{2}},d\bar{x}\rangle=0$ at every
point $P\in \Phi.$\medskip \newline(iv) \qquad$q_{1}=\lambda$\thinspace$q_{2}$,
where $\lambda$ is an arbitrary non-vanishing constant.
\end{proposition}

\end{document}